\documentclass[12pt]{article}
\usepackage{amssymb}
\def\qed{\hfill$\Box$}

\def\th{\theta}
\def\wt{\widetilde}
\def\F{\hbox{$I\hskip -4pt F$}}
\def\sF{\hbox{$\sc I\hskip -2.5pt F\!$}}
\def\E{\hbox{$I\hskip -4pt E$}}
\def\sE{\hbox{$\sc I\hskip -2.5pt E\!$}}
\def\DD{{\cal D}}
\def\Z{\hbox{$Z\hskip -5.2pt Z$}}
\def \N{\hbox{$I\hskip -4pt N$}}
\def\sN{\hbox{$\ssc I\hskip -2.5pt N$}}

\def\dis{\displaystyle}
\def\wh{\widehat}
\def\sc{\scriptstyle}
\def\ssc{\scriptscriptstyle}

\def\cl{\centerline}
\def\ol{\overline}
\def\ul{\underline}
\def\rar{\rightarrow}

\def\Rar{\Rightarrow}

\def\bs{\backslash}
\def\vs{\vspace}
\def\a{\alpha}
\def\b{\beta}
\def\d{\delta}
\def\e{\epsilon}
\def\g{\gamma}
\def\l{\lambda}
\def\p{\partial}
\def\G{\Gamma}
\def\stl{\stackrel}
\def\HH{{\cal H}}
\def\JJ{{\cal J}}
\def\AA{{\cal A}}
\def\BB{{\cal B}}
\def\WW{{\cal W}}
\def\II{{\cal I}}
\def\MM{{\cal M}}
\def\N{\mathbb{N}}
\def\sN{\mathbb{N}}
\def\Z{\mathbb{Z}}

\def\sF{\mathbb{F}}
\def\F{\mathbb{F}}
\def\E{\mathbb{E}}
\def\sE{\mathbb{E}}
\def\ni{\noindent}
\def\hi{\hangindent}
\def\ha{\hangafter}

\begin{document}
\ni {{\large\bf Classification of derivation-simple color algebras
related to locally finite derivations}\,\\[2pt] 
\cl{(Appeared in {\em J.~Math.~Phys.} {\bf45} (2004) 525--536.)}
\vs{-9pt}
\par \setcounter{page}{1} \ni\hi4.5ex\ha0 $\!\!\!\!\!\!$Yucai Su$^{\,a}$,
 \ Kaiming Zhao$^{\,b}$ \ and \ Linsheng
Zhu$^{\,c}$\vs{-2pt}\par\ni\hi4.5ex\ha0
 $\!\!\!\!\!\!${\it $^{a)}$ Department of Mathematics, Shanghai
Jiaotong University, Shanghai 200030,
China}\vs{-7pt}\par\ni\hi4.5ex\ha0
 $\!\!\!\!\!\!${\it
$^{b)}$ Department of Mathematics, Wilfrid Laurier University,
Waterloo, Ontario, Canada
 N2L 3C5, and Institute of Mathematics, Academy of Mathematics {\rm\&} system
Sciences, Chinese Academy of Sciences, Beijing 100080,
China}\vs{-7pt}\par\ni\hi4.5ex\ha0
 $\!\!\!\!\!\!${\it $^{c)}$ Department of Mathematics,
Changshu College, Jiangsu 215500, China} \vs{-3pt}\par\ \par
\ni\hi3ex\ha0 {\small We classify the pairs $(\AA,\DD)$ consisting
of an $(\e,\G)$-color-commutative associative algebra $\AA$ with
an identity element over an algebraically closed field $\F$ of
characteristic zero and a finite dimensional subspace $\DD$ of
$(\e,\G)$-color-commutative locally finite color-derivations of
$\AA$ such that $\AA$ is $\G$-graded $\DD$-simple and the
eigenspaces for elements of $\DD$ are $\G$-graded. Such pairs are
the important ingredients in constructing some simple Lie color
algebras which are in general not finitely-graded. As some
applications, using such pairs, we construct new explicit simple
Lie color algebras of generalized Witt type, Weyl type.}\par\ni
\vs{-5pt}\par\ni {\small\bf I. INTRODUCTION}\par
Lie color algebras, a notion first appeared in mathematical
physics,$^{1,3,5-7,15}$ are generalizations of Lie algebras and
Lie superalgebras. Let us start with the definition. Let $\F$ be
an algebraically closed field of characteristic zero and let $\G$
be an additive group. A {\it skew-symmetric bicharacter} of $\G$
is a map
 $\e:\G\times\G\rar \F^\times=\F\bs\{0\}$
satisfying
$$
\e(\l,\mu)=\e(\mu,\l)^{-1},\,\;
\e(\l,\mu+\nu)=\e(\l,\mu)\e(\l,\nu),\,\; \forall\,\l,\mu,\nu\in\G.
\eqno(1.1)$$
It is clear that
$$
\e(\l,0)=1,\;\;\forall\,\l\in\G.
\eqno(1.2)$$
Let
$L=\oplus_{\l\in\G}L_{\l}$ be a $\G$-graded $\F$-vector space. For
a nonzero homogeneous element $a$, denote by $\bar{a}$ the unique
group element in $\G$ such that $a\in L_{\bar{a}}$. We shall call
$\bar{a}$ the {\it color} of $a$. The $\F$-bilinear map
$[\cdot,\cdot]: L\times L\rar L$ is called a Lie color bracket on
$L$ if the following conditions are satisfied:
$$
\matrix{
 [a,b]=-\e(\bar{a},\bar{b})[b,a],\hfill&\mbox{(skew symmetry)}
\vs{4pt}\hfill\cr
 [a,[b,c]]=[[a,b],c]+\e(\bar{a},\bar{b})[b,[a,c]],\hfill&
\mbox{(Jacobi identity)}\cr}$$ for all homogeneous elements
$a,b,c\in L$. The algebra structure $(L,[\cdot,\cdot])$ is called
an {\it $(\e,\G)$-Lie color algebra} or simply a {\it Lie color
algebra}. If $\G=\Z/2\Z$ and $\e(i,j)=(-1)^{ij},\forall\,i,j\in
\Z/2\Z$, then $(\e,\G)$-Lie color algebras are simply Lie
superalgebras. For Lie color algebras, we refer the reader to the
book.$^1$
\par
For any $\G$-graded $\F$-vector space $V$, we denote
$$
H(V)=\{\mbox{ all homogeneous elements in }V\}.$$ Let
$\AA=\oplus_{\l\in\G}\AA_\l$ be a $\G$-graded associative
$\F$-algebra with an identity element 1, i.e.,
$\AA_\l\AA_\mu\subset\AA_{\l+\mu}$ for all $\l,\mu\in\G$. So $1\in
\AA_0$. We say that $\AA$ is {\it graded simple} if $\AA$ does not
have nontrivial $\G$-graded ideals. If we define the bilinear
product $[\cdot,\cdot]$ on $\AA$ by
$$
[x,y]=xy-\e(\bar{x},\bar{y})yx,\ \forall\,x,y\in H(\AA),
\eqno(1.3)$$ then $(\AA,[\cdot,\cdot])$ becomes a Lie color
algebra.
\par
A {\it Lie color ideal\,} $U$ of $\AA$ is a $\G$-graded vector space
$U$ of $\AA$ such that $[\AA,U]\subset U$. Sometimes it is called an
$(\e,\G)$-Lie ideal. The {\it $\e$-center} $Z_\e(\AA)$ of $\AA$ is
defined as
$$
Z_\e=Z_\e(\AA)=\{x\in\AA\,|\,[x,\AA]=0\}.$$ It is easy to see that
$Z_\e(\AA)$ is $\G$-graded. We say that $\AA$ is {\it
color-commutative} (or $\e$-color-commutative) if $Z_\e(\AA)=\AA$,
i.e., $[\AA,\AA]=0$.
\par
Let $\AA$ be an $(\e,\G)$-color-commutative associative algebra with an
identity element $1$. A nonzero
$\F$-linear transformation $\p:\AA\rar\AA$ is called a {\it
homogeneous color-derivation} of degree $\l\in\G$
if
$$
\matrix{
 \p(a)\in\AA_{\l+\mu},\ \forall\,a\in\AA_\mu,\ \mu\in\G\,\ {\rm and}
\vs{6pt}\hfill\cr
 \p(ab)=\p(a)b+\e(\l,\bar{a})a\p(b),\ \forall\,a,b\in H(\AA).
\hfill\cr} \eqno\matrix{(1.4)\!\vs{6pt}\cr\ \cr}$$ For
convenience, we shall often denote $\bar\p=\l$ if $\p$ has degree
$\l$. Clearly $\p(c)=0$ for all $c\in\F$. Denote ${\rm
Der}^\e(\AA)=\oplus_{\l\in\G}{\rm Der}_\l^\e(\AA)$, where ${\rm
Der}_\l^\e(\AA)$ is the $\F$-vector space spanned by all
homogeneous color derivations of degree $\l$. Similar to  the Lie
algebra case, it is easy to verify that ${\rm Der}_\l^\e(\AA)$
becomes a Lie color algebra under the Lie color bracket
$$
[\p,\p']=\p\p'-\e(\bar{\p},\bar{\p'})\p'\p,\ \forall\,\p,\p'\in
H({\rm Der}^\e(\AA)),$$ where $\p\p'$ is the composition of the
operators $\p$ and $\p'$.
\par
Let $\DD=\oplus_{\l\in\G}\DD_\l$ be an $(\e,\G)$-color-commutative
subspace of ${\rm Der}^\e(\AA)$, i.e.,
$$
\p\p'=\e(\bar{\p},\bar{\p'})\p'\p,\ \forall\,\p,\p'\in H(\DD).
\eqno(1.5)$$ Recall that the associative algebra $\AA$  is called
{\it graded $\DD$-simple} if $\AA$ has no nontrivial graded
$\DD$-stable ideals.$^4$
\par
A linear transformation $T$ on a vector space $V$ is called {\it
locally finite} if
$$
{\rm dim}({\rm span}\{T^m(v)\,|\, m\in\N\})<\infty,$$ for any
$v\in V$. The transformation $T$ is called {\it locally nilpotent}
if for any $v\in V$, we have $T^n(v)=0$ for some $n\in\N$, and $T$
is called {\it semi-simple} if it acts diagonalizably on $V$.
\par
For a pair $(\AA,\DD)$ of an $(\e,\G)$-color-commutative
associative algebra with an identity element and an
$(\e,\G)$-color-commutative subspace $\DD$ of ${\rm Der}^\e(\AA)$,
Passman$^{\,4}$ proved that the Lie color algebra (including the
Lie algebra case) $\AA\DD= \AA\otimes\DD$ is simple if and only if
$\AA$ is graded $\DD$-simple and $\AA\DD$ acts faithfully on $\AA$
(except a minor case). The authors of the present paper$^{\,11}$
(see also Refs.~9, 10 and 14) constructed (associative and Lie)
color algebras of Weyl type $\AA[\DD]$, which is the color
commutative algebra generated by $\AA$ and $\DD$ (as operators on
$\AA$), and proved that $\AA[\DD]$ is simple as an associative
algebra or is {\it central simple} as a Lie color algebra (i.e.,
the derived subalgebra modulo its $\e$-center is simple) if and
only if $\AA$ is graded $\DD$-simple (except a minor case in Lie
case). However, it is still a question of how to construct new
explicit simple Lie color algebras of generalized Witt type or
Weyl type.
\par
The problem of classifying all the pairs $(\AA,\DD)$ of a
commutative associative algebra $\AA$ with an identity element and
a finite-dimensional locally finite commutative derivation
subalgebra $\DD$ such that $\AA$ is {\it $\DD$-simple} (i.e.,
$\AA$ does not have $\DD$-stable ideals), was settled in Ref.~8
(using the pairs $(\AA,\DD)$, Xu constructed explicit simple Lie
algebras of generalized Cartan type$^{\,12}$ and of generalized
Block type$^{\,13}$). However, this problem becomes much more
complicated in color case.
\par
In order to construct explicit new simple Lie color algebras of
generalized Witt, Weyl types, the first aim of the present paper
is to give a classification of all the pairs $(\AA,\DD)$ of an
$(\e,\G)$-color-commutative associative algebra $\AA$ with an
identity element over an algebraically closed field $\F$ of
characteristic zero and a finite dimensional subspace $\DD$ of
$(\e,\G)$-color-commutative locally finite color-derivations of
$\AA$ such that $\AA$ is $\G$-graded $\DD$-simple and the
eigenspaces for elements of $\DD$ are $\G$-graded (see Theorem
2.2). Then in Section 3, as some applications, using the pairs
$(\AA,\DD)$, we construct explicit new simple Lie color algebras
(including Lie superalgebras) of generalized Witt, Weyl types (see
Theorem 3.1).
\par\ni
\vs{-5pt}\par\ni {\small\bf 2. $\DD$-SIMPLE COLOR ALGEBRAS}\par
In this section, we shall classify the pairs $(\AA,\DD)$ of an
$(\e,\G)$-commutating associative algebra $\AA$ with an identity
element $1$ and a finite-dimensional subspace $\DD$ of
$(\e,\G)$-commutative locally finite color derivations of $\AA$
such that $\AA$ is graded $\DD$-simple and the eigenspaces for
elements of $\DD$ are $\G$-graded. \par First we would like to
remark that the eigenspace of a derivation is not necessarily
$\G$-graded. Since we are considering $\G$-graded algebras, it is
natural that we require the eigenspaces for elements of $\DD$ are
$\G$-graded.\par We shall start with constructing explicitly such
pairs $(\AA,\DD)$. The motivation  to construct such pairs will
become clear in the proof of Theorem 2.2 below. Actaully, the
proof of Theorem 2.2 leads us the way to construct such pairs.
\par
Set
$$
\G_+=\{\l\in\G\,|\,\e(\l,\l)=1\},\;\;\;
\G_-=\{\l\in\G\,|\,\e(\l,\l)=-1\}. $$ Then by (1.1), $\G_+$ is a
subgroup of $\G$ with index $\le2$. For any graded subspace $\BB$
of $\AA$, we define $\BB_+=\oplus_{\l\in\G_+}\BB_\l$, then $\BB_+$
is $\G$-graded. Similarly we can define $\BB_-$. Since
$\G=\G_+\cup\G_-$, it follows that $\BB=\BB_+\oplus \BB_-$. By
(1.5), we have
$$
a^2=0\ \ \ \mbox{ or \ }\ \p^2=0\mbox{ \ if }\ \ol a\in\G_- \mbox{
\ or }\ \ol\p\in\G_-. \eqno(2.1)$$ For $m,n\in\Z$, we denote
$$
\ol{m,n}=\{m,m+1,\cdots,n\}.$$
\par
To construct the pair $(\AA,\DD)$, first we construct a {\it
$\G$-graded $\e$-commutative field extension} $\E$ of $\F$ (i.e.,
each nonzero homogeneous element of $\E$ is invertible). To do
this, let $\G^0\subset\G_+$ be a subgroup of $\G$ and let $\E_0$
be a field extension of $\F$. Let $e:\G^0\times\G^0\to\E_0^\times=
\E_0\bs\{0\}$ be a 2-variable function $e:(\a,\b)\mapsto
e_{\a,\b}$ such that
$$
e_{\a,\b}=\e(\a,\b)e_{\b,\a},\;\;e_{\a,0}=1,\;\;
e_{\a,\b}e_{\a+\b,\g}=e_{\a,\b+\g}e_{\b,\g},\;\;
\forall\,\a,\b,\g\in\G^0. \eqno(2.2)$$ You will see that these are
required by the associativity of the algebra we are going to
construct. Let $\E=\E_0[\G^0]={\rm
span}_{\sE_0}\{E_\a\,|\,\a\in\G^0\}$ be a $\G^0$-graded
$\e$-commutative associative algebra over $\E_0$ such that $E_\a$
has color $\ol E_\a=\a$, with the multiplication
$$
E_\a\cdot E_\b=e_{\a,\b}E_{\a+\b},\;\;\forall\,\a,\b\in\G^0.
\eqno(2.3)$$  From (2.2) it is easy to see that $\E$ is a
$\G$-graded $\e$-commutative field extension of $\F$.
\par
Let
$$
\ul k=(k_1,k_2,k_3,k_4)\in\N^4\mbox{ \ such that \
}k=k_1+k_2+k_3+k_4>0.$$ We also require that $k_4=0$ if
$\G_-=\emptyset$. We shall construct $\DD$ which will be spanned
by color derivations $\p_p,p\in\ol{1,k}$ such that
$$
\matrix{ \p_p\mbox{ is semi-simple with color }\ol\p_p=0,\hfill&
\forall\;p\in\ol{1,k_1}, \vs{6pt}\hfill\cr \p_{k_1+p}\mbox{ is
locally finite but not semi-simple with color }
\ol\p_{k_1+p}=0,\hfill& \forall\;p\in\ol{1,k_2}, \vs{6pt}\hfill\cr
\p_{k_1+k_2+p}\mbox{ is locally nilpotent with color }
\ol\p_{k_1+k_2+p}\in\G_+,\hfill& \forall\;p\in\ol{1,k_3},
\vs{6pt}\hfill\cr \p_{k_1+k_2+k_3+p}\mbox{ is locally nilpotent
with color } \ol\p_{k_1+k_2+k_3+p}\in\G_-,\hfill&
\forall\;p\in\ol{1,k_4}, \hfill\cr} \eqno \matrix{(2.4)\vs{6pt}\cr
(2.5)\vs{6pt}\cr(2.6)\vs{6pt}\cr(2.7)\cr}$$ (cf.~(2.21) and
(2.22)). To this end, we first need to construct $\AA$ which will
be the tensor product of two algebras $\AA=\AA_1\otimes\AA_2$
(cf.~(2.19)) such that $\AA_1$ is a ``group-algebra-like'' algebra
(cf.~(2.12)) and $\AA_2$ is a ``polynomial-like'' algebra
(cf.~(2.16)).
\par
Now we construct $\AA_1$ such that $\p_p|_{\AA_1}$ are nonzero
semi-simple operators for $p\in\ol{1,k_1+k_2}$ and
$\p_{k_1+k_2+p}|_{\AA_1}$ are zero operators for
$p\in\ol{1,k_3+k_4}$ (cf.~(2.4)-(2.7) and (2.21)-(2.22)). To do
this, let $G$ be a {\it nondegenerate} additive subgroup of
$\F^{k_1+k_2}$, i.e., $G$ contains an $\F$-basis of
$\F^{k_1+k_2}$. If $k_1+k_2=0$, we take $G=\{0\}$. An element in
$G$ is usually denoted by
$$
\ul a=(a_1,a_2,\cdots,a_k)\mbox{ with }a_p=0,\;\forall\,p>k_1+k_2.
\eqno(2.8)$$ Let $\wh{\ }:G\to\G_+$ be a map $\wh{\ }:\ul
a\mapsto\wh{\ul a}$ satisfying
$$
\wh 0=0,\;\;\; \th_{\ul a,\ul b}:= \wh{\ul a}+\wh{\ul b}
-\!\!\wh{\,\,\ul a\!+\!\ul b\,\,}\!\! \in\G^0,\;\;\forall\,\ul
a,\ul b\in G. \eqno(2.9)$$ Let $f(\cdot,\cdot):G\times
G\to\E_0^\times$ be a map such that
$$
f(\ul a,\ul b)=\e(\wh{\ul a},\wh{\ul b})f(\ul b,\ul a),\;\;\;
f(\ul a,0)=1, \eqno(2.10)$$
$$
e_{_{\th_{\ul a,\ul b},\th_{\ul a+\ul b,\ul c}}} f(\ul a,\ul
b)f(\ul a+\ul b,\ul c)= \e(\wh{\ul a},\th_{\ul b,\ul c})
e_{_{\th_{\ul b,\ul c},\th_{\ul a,\ul b+\ul c}}} f(\ul b,\ul
c)f(\ul a,\ul b+\ul c), \eqno(2.11)$$ for $\ul a,\ul b,\ul c\in
G$. Denote by $\AA_1=\AA(G,\E,f)$ the $(\e,\G)$-color commutative
associative algebra with $\E$-basis $\{x^{\ul a}\,|\,\ul a\in G\}$
or $\E_0$-basis $\{E_\a x^{\ul a}\,|\, (\a,\ul a)\in\G^0\times
G\}$ such that $x^{\ul a}$ has color $\wh{\ul a}$ and
$$
x^{\ul a}\cdot x^{\ul b}=f(\ul a,\ul b)E_{\th_{\ul a,\ul b}}x^{\ul
a+\ul b}, \;\;\forall\,\ul a,\ul b\in G, \eqno(2.12)$$ and in
general
$$
E_\a x^{\ul a}\cdot E_\b x^{\ul b}= \e(\wh{\ul
a},\b)e_{\a,\b}e_{\a+\b,\th_{\ul a,\ul b}} f(\ul a,\ul
b)E_{\a+\b+\th_{\ul a,\ul b}}x^{\ul a+\ul b},\;
\forall\,\a,\b\in\G^0,\,\ul a,\ul b\in G, \eqno(2.13)$$
(cf.~(2.3)). The $\e$-commutativity and associativity of $\AA_1$
are guaranteed by conditions (2.10) and (2.11).
\par
Now we shall construct $\AA_2$ such that $\p_p|_{\AA_2}=0$ for
$p\in\ol{1,k_1}$ and $\p_p|_{\AA_2}$ are nonzero locally nilpotent
operators for $p\in\ol{k_1+1,k}$ (cf.~(2.4)-(2.7) and
(2.21)-(2.22)). To this end, let $t_{k_1+1},\cdots,t_k$ be
$k_2+k_3+k_4$ variables such that each $t_p$ has color $\ol t_p$
satisfying
$$
\ol t_{k_1+p}=0,\;\;\;\ol t_{k_1+k_2+q}\in\G_+,\;\;\; \ol
t_{k_1+k_2+k_3+r}\in\G_-, \eqno(2.14)$$ for
$p\in\ol{1,k_2},\,q\in\ol{1,k_3},\,r\in\ol{1,k_4}$. For
convenience, we denote $t_p=0$ if $p\le k_1$. Denote
$\JJ=\{0\}^{k_1}\times\N^{k_2+k_3}\times\Z_2^{k_4}$, where
$\Z_2=\Z/2\Z$, i.e., $\JJ$ is the subset of $\F^k$ consisting of
the following elements:
$$
\ul i=(i_1,i_2,\cdots,i_k), \eqno(2.15)$$ with $i_p=0$ for $p\le
k_1$, and $i_q\in\N$ for $q\in\ol{k_1+1,k_1+k_2+k_3}$, and
$i_r=0,1$ for $q>k_1+k_2+k_3$, ((2.1) and (2.16) explain why we
shall have $i_q=0,1$ for $q>k_1+k_2+k_3$). Let
$\AA_2=\E[t_{k_1+1},\cdots,t_k]$ be the $\e$-commutative algebra
of polynomials in $k_2+k_3+k_4$ variables with an $\E$-basis
consisting of the elements
$$
t^{\ul i}=t_{k_1+1}^{i_{k_1+1}}\cdots t_k^{i_k},\;\;\forall\,\ul
i\in\JJ, \eqno(2.16)$$ or $\E_0$-basis $\{E_\a t^{\ul
i}\,|\,(\a,\ul i)\in\G^0\times\JJ\}$ such that
$$
E_\a t^{\ul i}\cdot E_\b t^{\ul j}=e_{\a,\b}
\prod_{p=k_1+1}^k\e(\ol t_p,\b)^{i_p} \prod_{k_1<p<q\le k}\e(\ol
t_q,\ol t_p)^{i_qj_p}E_{\a+\b}t^{\ul i+\ul j},
\;\forall\,\a,\b\in\G^0,\,\ul i,\ul j\in\JJ, \eqno(2.17)$$
(cf.~(2.3)), where we use the convention that $t^{\ul i}=0$ if
$\ul i\notin \JJ$. For convenience, we shall denote
$$
\e_{\ul i,\b}=\prod_{p=k_1+1}^k\e(\ol t_p,\b)^{i_p},\;\;\;
\wt\e_{\ul i,\ul j}=\prod_{k_1<p<q\le k}\e(\ol t_q,\ol
t_p)^{i_qj_p},\;\; \forall\,\ul i,\ul j\in\JJ,\,\b\in\G^0.
\eqno(2.18)$$
\par
{\it Definition 2.1.} We define $\AA=\AA(\ul k,G,\E,f)$ to be the
$(\e,\G)$-commutative associative algebra with the identity
element $1=E_0=x^0$, which is the tensor product of algebras
$\AA=\AA_1\otimes_{\sE}\AA_2$, having $\E_0$-basis
$$
E_\a x^{\ul a,\ul i}=E_\a x^{\ul a}t^{\ul i},\;\; \forall\,(\a,\ul
a,\ul i)\in\G^0\times G\times\JJ, \eqno(2.19)$$ with the
multiplication
$$
E_\a x^{\ul a,\ul i}\cdot E_\b x^{\ul b,\ul j}= \e_{\ul
i,\b}e_{\a,\b}\e_{\ul i,\wh{\ul b}}\wt\e_{\ul i,\ul j} \e(\wh{\ul
a},\b)e_{\a+\b,\th_{\ul a,\ul b}}f(\ul a,\ul b) E_{\a+\b+\th_{\ul
a,\ul b}}x^{\ul a+\ul b,\ul i+\ul j}, \eqno(2.20)$$ for
$\a,\b\in\G^0,\,\ul a,\,\ul b\in G,\,\ul i,\ul j\in\JJ$
(cf.~(2.3), (2.13), (2.17) and (2.18)).
\par
For $a\in\F^k,p\in\ol{1,k}$, we denote
$$
a_{[p]}=(0,\cdots,0,\stl{p}{a},0,\cdots,0)\in\F^k.$$ For
$p\in\ol{1,k}$, we define the linear transformations
$\p_p,\p_{t_p},\p_p^*$ on $\AA$ such that they have color $-\ol
t_p$ (in particular, they have color $0$ if $p\le k_1+k_2$,
cf.~(2.14)), and
$$
\p_p=\p^*_p+\p_{t_p}, \eqno(2.21)$$
$$
\p^*_p(E_\a x^{\ul a,\ul i})=a_pE_\a x^{\ul a,\ul i},\;\;\;
\p_{t_p}(E_\a x^{\ul a,\ul i}) =\e(\ol\p_{t_p},\a+\wh{\ul
a})\prod_{q=1}^{p-1}\e(\ol\p_{t_p},\ol t_q)^{i_q} {\sc\,}i_pE_\a
x^{\ul a,\ul i-1_{[p]}}, \eqno(2.22)$$ for $(\a,\ul a,\ul
i)\in\G^0\times G\times\JJ$. Clearly, $\p^*_p=0$ if $p>k_1+k_2$ by
(2.8), and $\p_q=0$ if $q\le k_1$ by (2.15). Then
$\p_p,\p^*_p,\p_{t_p}$ are $\e$-derivations of $\AA$ for
$p\in\ol{1,k}$. We call $\p^*_p$ a {\it grading operator} (or {\it
degree operator}), $\p_{t_p}$ a {\it down-grading operators}, and
$\p_p=\p^*_p+\p_{t_p}$ a {\it mixed operator} if both $p^*_p$ and
$\p_{t_p}$ are nonzero. Then
$$
\DD={\rm span}_{\sF}\{\p_p\,|\,p\in\ol{1,k}\}, \eqno(2.23)$$ is a
finite dimensional subspace of $\e$-commutative locally finite
color derivations of $\AA$ such that the eigenspaces for elements
of $\DD$ are $\G$-graded.
\par
{\bf Theorem 2.2}. {Let $\AA=\sum_{\a\in \G}\AA_\a$ be an
$\e$-commutative associative graded algebra with an identity
element over an algebraically closed field $\F$ of characteristic
zero and let $\DD=\sum_{\a\in \G}\DD_\a$ be a finite-dimensional
$\G$-graded $\F$-subspace of $\e$-commutative locally finite
color-derivations of $\AA$ such that the eigenspaces for elements
of $\DD$ are $\G$-graded. Then $\AA$ is graded $\DD$-simple if and
only if $\AA$ is isomorphic to the algebra of the form $\AA(\ul
k,G,\E,f)$ defined in (2.19) and (2.20), and $\DD$ is of the form
(2.21)-(2.23).}
\par
{\it Proof}. ``$\Leftarrow$'': Let $\II$ be a $\G$-graded
$\DD$-stable nonzero ideal of $\AA=\AA(\ul k,G,\E,f)$. By (2.21)
and (2.22), we see that
$$
\biggl({\sc\,}\bigcup_{(\a,\ul a)\in\G^0\times G}\F(E_\a x^{\ul
a}){\sc\,}\biggr)\bs\{0\},$$ is the set of the common eigenvectors
of $\DD$. We also see that if a homogeneous element $\p\in H(\DD)$
has a nonzero eigenvalue, then $\p\in\DD_0$. Thus
$\sum_{0\ne\a\in\G}\DD_\a$ acts locally nilpotently on $\II$.
Since $\DD_0$ is commutative (cf.~(1.2) and (1.5)), and $\DD_0$
commutes with $\sum_{0\ne\a\in\G}\DD_\a$, and
$\sum_{0\ne\a\in\G}\DD_\a$ is color-commutative, by linear
algebra, $\II$ must contain a common eigenvector of $\DD$. Thus
$E_\a x^{\ul a}\in\II$ for some $(\a,\ul a)\in\G^0\times G$. Then
$$
1=e^{-1}_{-\a,\a}f(-\ul a,\ul a)^{-1} (E_{-\a}x^{-\ul
a})\cdot(E_\a x^{\ul a})\in\II,$$ (cf.~(2.13)). Hence $\II=\AA$.
This proves that $\AA$ is graded $\DD$-simple.
\par
``$\Rightarrow$'': Suppose $\p\in H(\DD)$ has a nonzero eigenvalue
$a\in\F$ such that $u_a\in H(\AA)$ is a corresponding eigenvector.
Then we have $\p(u_a)=a u_a$, and so $\ol\p+\ol u_a=\ol u_a$ by
(1.4). Thus $\ol\p=0$. In other words, we have
$$
\p\in H(\DD),\,\ol\p\ne0\ \Rar\ \p\mbox{ \ acts locally nilpotent
on \ }\AA. \eqno(2.24)$$ Since $\F$ is algebraically closed and
$\DD$ is a finite dimensional subspace of $\e$-commutative locally
finite color derivations of $\AA$, from linear algebra, we have
$$
\AA=\bigoplus_{\ul a\in\DD^*}\AA(\ul a),$$ where $\DD^*$ is the
dual space of $\DD$, and
$$
\AA(\ul a)=\{u\in \AA\,|\,(\p-\ul a(\p))^m(u)=0\mbox{ for }\p\in
H(\DD)\mbox{ and some } m\in\N\},$$ for $\ul a\in\DD^*$ (note that
$\ul a(\p)=0$ if $\ol\p\ne0$ by (2.24)). Denote
$$
G=\{\ul a\in\DD^*\,|\,\AA(\ul a)\ne0\}.$$ By (2.24), $G$ can be
viewed as a subset of $\DD_0^*$ by the restriction $\ul
a\mapsto\ul a|_{\DD_0}$. For any $\ul a\in G,\,n\in\N$, we define
$$
\AA(\ul a)^{(n)}= \{u\in\AA\,|\,(d_1-\ul a(d_1))\cdots(d_{n+1}-\ul
a(d_{n+1}))(u)=0, \;\forall\,d_1,\cdots,d_{n+1}\in H(\DD)\}.
\eqno(2.25)$$ Then
$$
\AA(\ul a)=\bigcup_{n=0}^\infty \AA(\ul a)^{(n)},\;\forall\,\ul
a\in G.$$ A nonzero vector in $\AA(\ul a)^{(0)}$ is called a {\it
root vector} with root $\ul a$. For any homogeneous root vector
$u\in\AA(\ul a)^{(0)}$, clearly $\AA u$ is a $\G$-graded
$\DD$-stable ideal of $\AA$. Thus $\AA u=\AA$. In particular,
$vu=1$ for some $v\in\AA$. So any homogeneous root vector is
invertible. For a root vector $u\in H(\AA(\ul a)^{(0)})$ with $\ul
a\in G$ and any $\p\in H(\DD)$, we have
$$
\begin{array}{lll}
0&=&\p(1)=\p(uu^{-1})=\p(u)u^{-1}+\e(\ol\p,\ol u)u\p(u^{-1})
\vs{4pt}\\
&=&{\ul a}(\p)uu^{-1}+\e(\ol\p,\ol u)u\p(u^{-1})\vs{4pt}\\
&=&
\left\{\matrix{
\e(\ol\p,\ol u)u\p(u^{-1})\hfill&\mbox{if \ }\ol\p\ne 0,\vs{4pt}\hfill\cr
\ul a(\p)+u\p(u^{-1})\hfill&\mbox{if \ }\ol\p=0,\hfill\cr}\right.
\end{array}
$$ because $\ul a(\p)=0$ if $\ol\p\ne0$ by (2.24). This
implies
$$
\p(u^{-1})=-\ul a(\p)u^{-1}, \eqno(2.26)$$ by (2.24). Hence
$$
-\ul a\in G,\;\;\;\forall\;\ul a\in G. \eqno(2.27)$$ For any $x\in
H(\AA(\ul a)^{(0)}),\,y\in H(\AA(\ul b)^{(0)})$ and $\p\in
H(\DD)$, we have
$$
\p(xy)=\p(x)y+\e(\ol\p,\ol x)x\p(y)= \left\{\matrix{
0\hfill&\mbox{if \ }\ol\p\ne0,\vs{4pt}\hfill\cr (\ul a(\p)+\ul
b(\p))xy\hfill&\mbox{if \ }\ol\p=0. \hfill\cr }\right.
$$ Hence
$$
\AA(\ul a)^{(0)}\cdot\AA(\ul b)^{(0)}\subset\AA(\ul a+\ul
b)^{(0)}, \;\;\forall\;\ul a,\ul b\in G.$$ Considering the
invertibility of root vectors, we have
$$
\AA(\ul a)^{(0)}\cdot\AA(\ul b)^{(0)}=\AA(\ul a+\ul b)^{(0)},
\;\;\;\forall\;\ul a,\ul b\in G.$$ In particular, we obtain
$$
\ul a+\ul b\in G,\;\;\;\forall\;\ul a,\ul b\in G. \eqno(2.28)$$
Thus by (2.27) and (2.28), $G$ is an additive subgroup of $D^*$.
Set
$$
\E=\AA(0)^{(0)}. \eqno(2.29)$$ Then $\E$ is a $\G$-graded field
extension of $\F$ such that $\E_0$ is a field extension of $\F$.
We set
$$
\G^0=\{\a\in\G\,|\,\E_\a\ne\{0\}\}.$$ Clearly, $\G^0$ is a
subgroup of $\G$ and $\G^0\subset\G_+$ by (2.1). For any
$\a\in\G^0$, choose $E_\a=1$ if $\a=0$, and $E_\a\in\E_\a\bs\{0\}$
if $\a\ne0$. Then $\{E_\a\,|\,\a\in\G^0\}$ forms an $\E_0$-basis
of $\E$. Thus we have (2.3) such that the coefficient $e_{\a,\b}$
satisfies (2.2) by color commutativity and associativity.
\par
First assume that $\AA(0)\ne\E$. Since ${\ul a}(\p)=0$ for any
homogeneous derivation $\p$ with $\ol\p\ne0$, for $u\in
\AA(\ul a)^{(m)},\,v\in\AA(\ul b)^{(n)}$ and
$d_1,\cdots,d_{m+n+1}\in H(\DD)$, by induction on $m+n+1$, we can write
$$
(d_1-(\ul a+\ul b)(d_1))\cdots(d_{m+n+1}-(\ul a+\ul
b)(d_{m+n+1}))(uv), \eqno(2.30)$$ as a linear combination of the
forms
$$
(d_{i_1}-\ul a(d_{i_1}))\cdots(d_{i_r}-\ul a(d_{i_r}))(u)\cdot
(d_{j_1}-\ul a(d_{j_1}))\cdots(d_{j_s}-\ul a(d_{j_s}))(v),
\eqno(2.31)$$ where
$$
r+s=m+n+1,\;\;\{i_1,\cdots,i_r,j_1,\cdots,j_s\}=\{1,\cdots,m+n+1\}.
$$
By definition (2.25), we obtain that (2.31) is zero, and so is
(2.30). It follows that $uv\in\AA(\ul a+\ul b)^{(m+n)}$. Thus
$$
\AA(\ul a)^{(m)}\cdot\AA(\ul b)^{(n)}\subset \AA(\ul a+\ul
b)^{(m+n)},\;\;\;\forall\;\ul a,\ul b\in G,\;m,n\in\N.
\eqno(2.32)$$ In particular, (since homogeneous root vectors are
invertible),
$$
\E\AA(\ul a)^{(m)}=\AA(\ul a)^{(m)}=\AA(\ul
a)^{(0)}\AA(0)^{(m)},\;\;\;\forall\;\ul a\in G,\; m\in\N,
\eqno(2.33)$$ (cf.~(2.29)). Hence each $\AA(\ul a)^{(m)}$ is a
vector space over the graded field $\E$. For any
$v\in\AA(0)^{(1)}$, we have $\DD(v)\subset\E$ and
$$
\DD(v)=0\ \Leftrightarrow\ v\in\E. \eqno(2.34)$$ Set
$$
\HH=\E\DD,\;\;\;\HH_1=\{\p\in\HH\,|\,\p(\AA(0)^{(1)})=\{0\}\},\;\;\;
k_1={\rm dim}_{\sE}\HH_1. \eqno(2.35)$$ Expression (2.34) implies
that $\AA(0)^{(1)}/\E$ is isomorphic to a subspace of the space
${\rm Hom}_{\sE}(\HH,\E)$ over $\E$. By linear algebra, there
exist subsets
$$
\{\p_{k_1+1},\p_{k_1+2},\cdots,\p_k\}\subset H(\DD),\;\;\;
\{t_{k_1+1},t_{k_1+2},\cdots,t_k\}\subset H(\AA(0)^{(1)}),
\eqno(2.36)$$ for some $k\in\N$, such that
$$
\AA(0)^{(1)}=\E+\sum_{l=k_1+1}^k\E t_l,\;\;\; \p_p(t_q)=\d_{p,q},
\;\;\;\forall\;p,q\in\ol{k_1+1,k}. \eqno(2.37)$$ Set
$$
\HH_2=\sum_{p=k_1+1}^k \E\p_p.$$ Then we have
$$
\HH=\HH_1\oplus\HH_2.$$ For convenience, denote
$$
t^{\ul i}=t_{k_1+1}^{i_{k_1+1}}\cdots t_k^{i_k}\mbox{ \ for \ }
\ul i=(i_{k_1+1},\cdots,i_k)\in\N^\ell,$$ where $\ell=k-k_1$.  By
(2.1) then
$$
t^{\ul i}=0\mbox{ \ if \ }i_p\ge2\mbox{ \ with \ } \ol
t_p\in\G_-\mbox{ \ for some \ } p\in\ol{k_1+1,k},$$ and
$$
t^{\ul i}\cdot t^{\ul j}= \prod_{k_1+1\le p<q\le k}\e(\ol t_q,\ol
t_p)^{i_qj_p} t^{\ul i+\ul j},\;\;\;\forall\;\ul i,\ul
j\in\N^\ell.$$ Furthermore, by (2.37), we can deduce by induction
on the level $|\ul i|:=\sum_{p=k_1+1}^k i_p$ that
$$
t^{\ul i}\in\AA(0)^{({\ssc\,}|\ul i|{\ssc\,})},\;\;\;\;
\p_p(t^{\ul i})=\prod_{k_1<q<p}\e(\ol
p_p,t_q)^{i_q}{\sc\,}i_pt^{\ul i-1_{[p]}}, \eqno(2.38)$$ for $\ul
i\in\N^\ell,p\in\ol{k_1+1,k}$. Set
$$
\wt\AA(0)=\sum_{\ul i\in\sN^\ell}\E t^{\ul i}\subset\AA(0).
\eqno(2.39)$$ Then $\wt\AA(0)$ forms a subalgebra of $\AA$. We
want to prove that $\AA(0)=\wt\AA(0)$. By (2.37),
$\AA(0)^{(1)}\subset\wt\AA(0)$. Suppose
$\AA(0)^{(m)}\subset\wt\AA(0)$ for some $1\le m\in\N$. By (2.25),
$\p(\AA(0)^{(m+1)})\subset\AA(0)^{(m)}\subset\wt\AA(0)$ for any
$\p\in\HH$. Thus, for $\p_1\in H(\HH)$ and $u\in H(\AA(0)^{(m)})$,
we may assume that
$$
\p_{k_1+1}(u)=\sum_{\ul i\in\sN^\ell} c_{\ul i}t^{\ul i},
\eqno(2.40)$$ where $c_{\ul i}\in H(\E)$ and $c_{\ul i}=0$ for all
but a finite number of $\ul i$. If $\ol\p_{k_1+1}\in\G_+$, then we
set
$$
u_1= \sum_{\ul i\in\sN^\ell}c_{\ul i}\e(\ol\p_{k_1+1},\ol c_{\ul
i})^{-1} (i_{k_1+1}+1)^{-1}t^{\ul i+1_{[k_1+1]}}\in H(\wt\AA(0)),
\eqno(2.41)$$ and we obtain
$$
\p_{k_1+1}(u)=\p_{k_1+1}(u_1). \eqno(2.42)$$ If
$\ol\p_{k_1+1}\in\G_-$, then by (2.1), $\p_{k_1+1}^2=0$, we must
have
$$
i_{k_1+1}=0\mbox{ \ \ if \ \ }c_{\ul i}\ne0, \eqno(2.43)$$
otherwise if (2.43) does not hold, then by (2.38) and (2.40) we
would have $\p_{k_1+1}^2(u)\ne0$, leading to a contradiction to
the fact  that $\p_{k_1+1}^2=0$. Thus we can still choose $u_1$ as
in (2.41) to give (2.42). Similarly, since $\p_{k_1+2}(u-u_1)\in
\AA(0)^{(m)}\subset\wt\AA(0)$, there exists $u_2\in H(\wt\AA(0))$
such that
$$
\p_{k_1+2}(u-u_1)=\p_{k_1+2}(u_2). \eqno(2.44)$$ Assume that
$u_2=\sum_{\ul i\in\sN^\ell}c'_{\ul i}t^{\ul i}$, where $c'_{\ul
i}\in H(\E)$. Since $\HH$ is color commutative, by (2.42) and
(2.44), we have
$$
\matrix{ 0\!\!\!\!&=\p_{k_1+1}\p_{k_1+2}(u_2)
\vs{4pt}\hfill\cr&\dis =\sum_{\ul i\in\sN^\ell}c'_{\ul i}
\e(\ol\p_{k_1+1}+\ol\p_{k_1+2},\ol c'_{\ul i})
e(\ol\p_{k_1+2},t_{k_1+1})^{i_{k_1+1}}
{\ssc\,}i_{k_1+1}i_{k_1+2}t^{\ul i-1_{[k_1+1]}-1_{[k_1+2]}}.
\hfill\cr}$$ Thus $i_{k_1+1}i_{k_1+2}=0$ if $c'_{\ul i}\ne0$.
Hence we can re-choose $u_2\in H(\wt\AA(0))$ such that
$$
\p_{k_1+1}(u_2)=0,\;\;\;\;\p_{k_1+2}(u-u_1)=\p_{k_1+2}(u_2).
$$ Similarly, we can find $u_2,\cdots,u_\ell\in
H(\wt\AA(0))$ such that
$$
\p_{k_1+p}(u-\sum_{q=1}^p u_q)=0,\;\;\;\;
\p_{k_1+1}(u_p)=\p_{k_1+2}(u_p)=\cdots=\p_{k_1+p-1}(u_p)=0,\;\;\;
\forall\;p\in\ol{2,\ell},$$ by induction on $p$. Thus we have
$$
\p_{k_1+p}(u-\sum_{q=1}^\ell u_q)=0,\;\;\;\;
\forall\;p\in\ol{1,\ell}. \eqno(2.45)$$ For any $\p,\p'\in
H(\HH_1)$, using (2.35) and (2.39) we deduce
$$
\p\p' (u-\sum_{p=1}^\ell u_p)\in\p(\AA(0)^{(m)})+\p'(\AA(0)^{(m)})
\subset\p(\wt\AA(0))+\p'(\wt\AA(0))=\{0\}. \eqno(2.46)$$ Now
(2.45) and (2.46) show that $u-\sum_{p=1}^\ell
u_p\in\AA(0)^{(1)}$. Thus by (2.35),
$$
\p(u-\sum_{p=1}^\ell u_p)=0,\;\;\;\;\forall\;\p\in H(\HH_1).
\eqno(2.47)$$ Then (2.45), (2.47) and the definition (2.25) show
that
$$
u-\sum_{p=1}^\ell u_p\in\AA(0)^{(0)}=\E.$$ Thus $u\in\wt\AA(0)$.
This proves
$$
\AA(0)=\wt\AA(0).$$ The case $\AA(0)=\E$ can be viewed as in the
general case $\AA(0)=\wt\AA(0)$ with $\ell=0$.
\par
We re-choose $\p_p,t_p,\,p\in\ol{1,k}$ as follows: Choose a
homogeneous $\F$-basis $\{\p_1,\cdots,\p_{k_1}\}$ of
$\DD\cap\HH_1$, and set $t_p=0$ for $p\in\ol{1,k_1}$, then $\p_p$
are semi-simple derivations on $\AA$ by (2.33) and (2.35). Let
$\ell_1$ be the dimension of the maximal locally nilpotent
$\F$-subspace of $\DD$. Clearly $\ell_1\le\ell=k-k_1$. Let
$k_2=\ell-\ell_1$. Now we choose $\p_{k_1+k_2+1},\cdots, \p_k$ to
be homogeneous locally nilpotent derivations of $\DD$ such that
the first $k_3$ derivations have colors in $\G_+$ and the last
$k_4$ derivations have colors in $\G_-$ for some $k_3,k_4$ with
$k_3+k_4=\ell_1$. Extend
$\{\p_p\,|\,p\in\ol{1,k_1}\cup\ol{k_1+k_2+1,k}\}$ to a homogeneous
$\F$-basis $\{\p_p\,|\,p\in\ol{1,k}\}$ of $\DD$. By the choices of
$\p_p$, then there exists $t_p\in\AA(0)^{(1)}$ for each
$p\in\ol{k_1+1,k}$ satisfying (2.36) and (2.37).
\par
For any $\ul a\in G$, we identify
$$
\ul a\ \leftrightarrow\ (\ul a(\p_1),\cdots,\ul a(\p_{k_1+k_2}))
\in\F^{k_1+k_2}.$$ Then $G$ is a nondegenerate subgroup of
$\F^{k_1+k_2}$ (otherwise, there exists
$\p\in\sum_{p=1}^{k_1+k_2}\F\p_p$ such that $\ul a(\p)=0$ for all
$\ul a\in G$ and so $\p$ is locally nilpotent, which contradicts
the maximality of $\ell_1$). Taking homogeneous root vector
$u\in\AA^{(0)}(\ul a)$, by (2.26) and (2.32), we have
$$
u^{-1}\AA(\ul a)\subset\AA(0),\;\;u\AA(0)\subset\AA(\ul a).
$$ Hence
$$
u\AA(0)=\AA(\ul a). \eqno(2.48)$$ In particular,
$$
\AA(\ul a)^{(0)}=\E u, \eqno(2.49)$$ is one-dimensional over $\E$.
Choose
$$
x^0=1,\;\;0\ne x^{\ul a}\in\AA(\ul a)^{(0)}
\;\;\;\mbox{for}\;\;0\ne\ul a\in G,$$ such that $x^{\ul a}$ is
homogeneous with color denoted by $\wh{\ul a}$. Since $x^{\ul a}$
is invertible, we have $\wh{\ul a}\in\G_+$. Then we have a map
$\wh{\ }$ satisfying (2.9). By (2.32) and (2.49), we have (2.12)
with $f(\ul a,\ul b)$ satisfying (2.10) and (2.11) by color
commutativity and associativity. By (2.48), we obtain
$$
\AA=\AA^{(0)}\AA(0)\cong\AA^{(0)}\otimes\AA(0),$$ where
$\AA^{(0)}=\oplus_{\ul a\in G}\AA(\ul a)^{(0)}$ is isomorphic to
the algebra $\AA_1$ defined in (2.12) and (2.13), and
$\AA(0)=\wt\AA(0)$ is isomorphic to the algebra $\AA_2$ defined in
(2.16) and (2.17). Therefore, the algebra $\AA$ is isomorphic to
the algebra $\AA(\ul k,G,\E,f)$ defined in (2.19) and (2.20), and
$\DD$ is of the form (2.23). This completes the proof of Theorem
2.2. \qed\par

\par\ni
\vs{-3pt}\par\ni {\small\bf 3. CONSTRUCTING SIMPLE LIE COLOR
ALGEBRAS FROM $\DD$-SIMPLE COLOR ALGEBRAS}\par
In this section, as applications, we shall construct some explicit
simple Lie color algebras using the pairs $(\AA,\DD)$ given in the
last section. For simplicity, we assume that the pairs $(\AA,\DD)$
in (2.19) and (2.23) satisfies
$$
\{u\in\AA\,|\,\DD(u)=0\}=\F.$$ This is equivalent to that
$\E_0=\F$ and $\G^0=\{0\}$.  So the map $\wh{\ }:G\to\G_+$ in
(2.9) is a group homomorphism and $\th_{\ul a,\ul b}=0$ for all
$\ul a,\ul b\in G$. In this case, noting that $\F$ is
algebraically closed, we prove that we can choose suitable basis
$\{x^{\ul a}\,|\,\ul a\in G\}$ such that the coefficient $f(\ul
a,\ul b)$ determined by (2.12), which satisfies (2.10) and (2.11),
has the following form:
$$
f(\ul a,\ul b)=\e(\wh{\ul a},\wh{\ul b})^{1\over2},\;\;\;\forall\;
\ul a,\ul b\in G, \eqno(3.1)$$ where the right-hand side is a
fixed square root such that (2.10) and (2.11) hold.
\par
Let $G'$ be a maximal subgroup of $G$ such that $x^{\ul a},\,\ul
a\in G'$ can be chosen so that (3.1) holds for $\ul a,\ul b\in
G'$. Suppose $G'\ne G$. Let $\ul c\in G\bs G'$ and set
$G''=G'+\Z\ul c$. If $G'\cap\Z\ul c=\{0\}$, we choose any $x^{\ul
c}\ne0$, and set $x^{\ul a+k\ul c}=\e(\wh{\ul a},\wh{\ul
c})^{-{k\over2}}x^{\ul a}\cdot(x^{\ul c})^k$ for $\ul a+k\ul c\in
G''$. If $G'\cap\Z\ne\{0\}$, then $G'\cap\Z\ul c=\Z\ul d$ for some
$\ul d= m\ul c,\,m>1$. In this case, since $\F$ is algebraically
closed, we can choose $x^{\ul c}$ such that $(x^{\ul c})^m=x^{\ul
d}$, and set $x^{\ul a+k\ul c}$ as above. In any case, the
coefficient $f(\ul a,\ul b)$ determined by (2.12) satisfies (3.1)
for $\ul a,\ul b\in G''$. But $G'\ne G''\supset G'$. This
contradicts the maximality of $G'$. This proves (3.1).
\par
Let $\F[\DD]$ be the $(\e,\G)$-commutative associative algebra with basis
$$
\{\p^\mu=\p_1^{\mu_1}\cdots\p_k^{\mu_k}\,|\,
\mu=(\mu_1,\cdots,\mu_k)\in \MM\},$$ where
$\MM=\N^{k_1+k_2+k_3}\times\Z_2^{k_4}$. For convenience, we denote
$\p^\mu=0$ if $\mu\notin\MM$. Denote
$$\matrix{
W=W(\ul k,G)=\AA\otimes\DD={\rm span}\{x^{\ul a,\ul i}\p_p\,|\,
(\ul a,\ul i)\in G\times\JJ,p\in\ol{1,k}\},\vs{6pt}\cr \WW=\WW(\ul
k,G)=\AA\otimes\F[\DD]={\rm span}\{x^{\ul a,\ul i}\p^\mu\,|\, (\ul
a,\ul i,\mu)\in G\times\JJ\times\MM\}.}
$$ Then as spaces, we have $W\subset\WW$. By regarding $\WW$ as operators on $\AA$, $\WW$
becomes a $\G$-graded associative algebra whose multiplication is
the composition of operators. Thus $\WW$ forms an $(\e,\G)$-Lie
color algebra under the bracket (1.3). We call $\WW$ a {\it Lie
color algebra of (generalized) Weyl type}. Clearly $\F$ is the
center of $\WW$. Let $\wt{\WW}=\WW/\F$ and let
$\ol\WW=[\wt\WW,\wt\WW]$ the derived algebra of $\wt\WW$.
Obviously, $W$ forms an $(\e,\G)$-Lie color subalgebra of $\WW$,
called a {\it Lie color algebra of (generalized) Witt type}. Using
results in Refs.~4 and 11, we obtain
\par
{\bf Theorem 3.1}. {The Lie color algebras $\ol{\WW}$ and $W$ are
simple if $k_1+k_2+k_3>0$ or $k_4>1$. Furthermore, $\ol\WW=\wt\WW$
if $k_1+k_2+k_3>0$ or otherwise, $\wt\WW=\ol\WW+\F t^{\ul
n}\p^\l$, where $\ul n$ and $\l$ are the largest elements
respectively in $\JJ$ and in $\MM$. } \qed\par Note that in case
$k=k_4=1$, $\ol{\WW}=0$ and $W=\F t_1\p_1$ are not simple. If
$k=k_4>1$, then we obtain finite dimensional simple Lie color
algebras $\ol\WW$ and $W$ of dimensions $2^{2n}-2$ and $n2^n$. In
particular, if $\G=\Z_2,\e(i,j)=(-1)^{ij},i,j\in\Z_2$, we obtain
the finite dimensional simple Lie superalgebras $\ol\WW=H(2n)$ and
$W=W(n)$ (see Ref.~2).
\def\EXAMPLE{\par
Next we construct simple Lie color algebras of (generalized) Block type.
To do this, we take $k=2$ in the construction of $\AA$, which is now
denoted by $\AA_2$.
To give explicit construction, we give new settings below.
Pick $\JJ_p\in\{\{0\},\Z_2,\N\}$ for $p=1,2$,
and set $\JJ=\JJ_1\times\JJ_2$.
Take an additive subgroup $G$ of $\F^2$ such that
$$
G_p\ne\{0\}\mbox{ if }\JJ_p=\{0\},\mbox{ \ and \ }
G_p=\{0\}\mbox{ if }\JJ_p=\Z_2
\ \ \mbox{for}\ \ p=1,2,
\eqno(3.6)$$
where $G_p=\{a_p\,|$ some $(a_1,a_2)\in G\}$ for $p=1,2$. Furthermore,
we require that
$$
\JJ_1\ne\JJ_2\mbox{ if }
G_1\ne\{0\}\ne G_2\mbox{ but }G\subset\F^2\mbox{ is degenerate}.
\eqno(3.7)$$
\par
For $p=1,2$, if $\JJ_p\bs\{0\}$,
we let $t_p$ be a variable with color $\ol t_p\in\G$ such that
$$
\left\{\matrix{
\ol t_p=0\hfill&\mbox{if \ }G_p\ne\{0\},\vs{4pt}\hfill\cr
\ol t_p\in\G_-\hfill&\mbox{if \ }\JJ_p=\Z_2,\vs{4pt}\hfill\cr
\ol t_p\in\G_+\hfill&\mbox{if \ }\JJ_p\ne\Z_2.
}\right.
\eqno(3.8)$$
We set $t_p=0$ if $t_p$ is not defined.
\par
Let $\wh{\ }$ and $f(\ul a,\ul b)$ satisfy (3.1). Let
$\AA_2=\AA_2(G,\JJ)$ be the $(\e,\G)$-commutative associative
algebra with a basis $\{x^{\ul a,\ul i}=x^{\ul
a}t_1^{i_1}t_2^{i_2} \,|\,(\ul a,\ul i)\in G\times\JJ\}$ and the
algebraic operation $\cdot$ defined by:
$$
x^{\ul a,\ul i}\cdot x^{\ul b,\ul j}=
f(\ul a,\ul b)\e(i_1\ol t_1+i_2\ol t_2,\wh{\ul b})
\e(i_2\ol t_2,j_1\ol t_1)x^{\ul a+\ul b,\ul i+\ul j},
\;\;\;\;\forall\;(\ul a,\ul i),(\ul b,\ul j)\in G\times\JJ,
\eqno(3.9)$$
such that $x^{\ul a,\ul i}$ has color
$\ol{x^{\ul a,\ul i}}=
\wh{\ul a}+i_1\ol t_1+i_2\ol t_2$.
We define the
derivations $\p_p$ of $\AA_2$ such that $\p_p$ has color $-\ol t_p$ and
$$
\p_p(x^{\ul a,\ul i})=a_p x^{\ul a,\ul i}+
e i_p x^{\ul a,\ul i-1_{[p]}},\;\;\;\;\forall\;
(\ul a,\ul i)\in G\times\JJ,\ p=1,2,
\eqno(3.10)$$
where $e=\e(\ol\p_1,\wh{\ul a})$
if $p=1$ and $e=\e(\ol\p_2,\wh{\ul a}+i_1\ol t_1)$ if $p=2$, and where
$x^{\ul a,\ul i}=0$ if $\ul i\notin\JJ$.
\par
To obtain the Lie color algebra structure
$(\AA_2,[\cdot,\cdot])$, we need to re-define colors
for elements in the Lie color algebra $(\AA_2,[\cdot,\cdot])$.
To distinguish colors between $(\AA_2,\cdot)$ and $(\AA_2,[\cdot,\cdot])$,
we use $\wt u$ to denote the color of $u$ in $(\AA_2,[\cdot,\cdot])$. We
define
$$
\wt u=\ol u+\ol\p_1+\ol\p_2,\;\;\;\;\forall\;u\in H(\AA_2).
\eqno(3.11)$$
We define the following algebraic operation $[\cdot,\cdot]$
on $\AA_2$:
$$
[u,v]=\p_1(u)\p_2(v)-\e(\wt u,\wt v)\p_1(v)\p_2(u),
\;\;\;\forall\;u,v\in H(\AA_2).
\eqno(3.14)$$
Note that we have $\wt{\p_1(u)\p_2(v)}=\wt u\wt v$ by (3.11), thus (3.14)
defines a $\G$-graded algebra structure on $\AA_2$.
It is straightforward to verify that $(\AA_2,[\cdot,\cdot])$ forms a
Lie color algebra such that $\F$ is the center of $(\AA_2,[\cdot,\cdot])$.
Denote $\BB_2=\BB_2(G,\JJ)=\AA/\F$.
\par
{\bf Theorem 3.2}. The Lie color algebra $\BB_2$ is simple if
${\rm dim\ssc\,}\BB_2\ge4$, i.e., $\JJ\ne\Z_2\times\Z_2$.
\par
{\it Proof of Theorem 3.2}.
We shall still use $x^{\ul a,\ul i}$ to denote elements in $\BB_2$.
The theorem can be proof in the following 3 cases:
(1) $\ol\p_1=\ol\p_2=0$ (2) $\ol\p_1=0,\ol\p_2\ne0$,
(3) $\ol\p_1\ne0\ne\ol\p_2.$ The proof of case (1) is exactly analogous to
that of Lie case [?]. We leave the other case to the reader.
\qed\par
As an example, suppose $G=\{0\},\JJ=\N\times\Z_2$.
}
\par
Using the pair $(\AA,\DD)$, one might construct other simple Lie
color algebras, for example, other series of Lie color algebras of
Cartan type.
 \small\vs{-7pt}\par\
\par\ni{\bf ACKNOWLEDGEMENTS}\vs{-1pt}\par YS and LZ are partially supported
by Morningside Center of Mathematics, Academy of Mathematics and
Systems Sciences, Chinese Academy of Sciences, Beijing 100080,
P.~R.~China during their visits to the Center. They wish to thank
the Center for the hospitality and the support. This work was
supported by NSF grant 10171064 of China,
two grants ``Excellent Young Teacher Program'' and
``Trans-Century Training Programme Foundation for the Talents''
from Ministry of Education of China
and the Foundation of Jiangsu Educational Committee.
\vs{-7pt}\par\
\par\ni{\bf REFERENCES}\vs{-2pt}\par
\def\r{\vs{-4pt}\
\par\ni\hi2.5ex\ha1}
\ni\hi2.5ex\ha1
\small
$^1\,$Yu. Bahturin, A. Mikhaliv, V. Petrogradskii and M. Zaicev,
  {\it Infinite dimensional Lie superalgebras}, Expos.~Math.
    Vol.~{\bf7}, de Gruyter, Berlin, 1992.
\r$^2\,$V.~G.~Kac, ``Lie superalgebras,'' Adv.~ Math. {\bf26}
   (1977), 8-96.
\r$^3\,$J.~Lukierski, V.~Rittenberg, ``Color-de Sitter and
  color-conformal superalgebras,'' Phys.~Rev.~D (3) {\bf18} (1978),
  385-389.
\r$^4\,$D.~P.~Passman, ``Simple Lie color algebras of Witt type,''
   J.~Algebra {\bf 208} (1998), 698-721.
\r$^5\,$V.~Rittenberg, D.~Wyler, ``Sequences of $Z_2\oplus Z_2$
  graded Lie algebras and superalgebras,'' J.~Math. Phys. {\bf19}
  (1978), 2193-2200.\r$^6\,$V.~Rittenberg, D.~Wyler, ``Generalized
  superalgebras,'' Nuclear Phys.~B {\bf139} (1978),
  189-202.
\r$^7\,$M.~Scheunert, R.~B.~Zhang, ``Cohomology of Lie
  superalgebras and their generalizations,'' J.~Math. Phys. {\bf39}
  (1998), 5024-5061.
\r$^8\,$Y.~Su, X.~Xu and H.~Zhang,
   ``Derivation-simple algebras and the structures of Lie algebras
   of Witt type,'' J.~Algebra {\bf 233} (2000), 642-662.
\r$^9\,$Y.~Su, K.~Zhao, ``Simple algebras of Weyl type,''
    Science in China A {\bf30} (2000) 1057-1063.
\r$^{10\,}$Y.~Su, K.~Zhao, ``Isomorphism classes and automorphism
  groups
    of algebras of Weyl type,'' Science in China A {\bf 45} (2002), 953-963.
\r$^{11\,}$Y.~Su, K.~Zhao and L.~Zhu,
    ``Simple Lie color algebras of Weyl types,'' Israel J. Math.,
    in press.
\r$^{12\,}$X.~Xu, ``New generalized simple Lie algebras of Cartan
    type over a field with characteristic 0,''
   J.~Algebra {\bf 224} (2000), 23-58.
\r$^{13\,}$X.~Xu, ``Generalizations of Block algebras,''
   Manuscripta Math. {\bf100} (1999), 489-518.
   \r$^{14}\,$K.~Zhao,
    ``Simple algebras of Weyl type II,"
   Proc. Amer. Math. Soc. {\bf130} (2002), 1323-1332.
\r$^{15\,}$A.~Zinoun, J.~Cortois, ``Pertinent statistics for
  describing a system of particles,'' J.~Math.~Phys. {\bf32} (1991),
  247-249.

\end{document}